\documentclass[11pt,reqno]{amsart}
\usepackage[utf8]{inputenc}
\usepackage{enumerate}
\usepackage{amssymb}
\usepackage{amsmath}
\usepackage{mathrsfs}
\usepackage{amsthm}
\usepackage{mathpazo}
\usepackage{extarrows}
\usepackage{color}
\usepackage{setspace}
\usepackage[colorlinks,linkcolor = blue,anchorcolor = red,urlcolor = blue,citecolor = blue]{hyperref}
\usepackage{graphicx}
\usepackage[square, comma, sort&compress, numbers]{natbib}

\theoremstyle{definition}
\newtheorem{Def}{Definition}[section]

\newtheorem{Rem}[Def]{Remark}

\theoremstyle{plain}
\newtheorem{Thm}[Def]{Theorem}
\newtheorem{Lem}[Def]{Lemma}
\newtheorem{Pro}[Def]{Proposition}
\newtheorem{Cor}[Def]{Corollary}

\newtheorem*{Que}{Question}
\usepackage[left=2.6cm, right=2.6cm, top=3cm, bottom=3cm]{geometry}
\usepackage{bm}
\usepackage{comment}

\def\IN{\mathbb N}\def\IR{\mathbb R}\def\IC{\mathbb C}\def\IZ{\mathbb Z}
\def\B{\mathcal B}\def\C{\mathcal C}\def\F{\mathcal F}\def\P{\mathcal P}\def\Q{\mathcal Q}\def\U{\mathcal U}\def\V{\mathcal V}\def\H{\mathcal H}

\def\asdim{\operatorname{asdim}}
\def\dimnuc{\operatorname{dim}_{\textup{nuc}}}

\def\diag{\operatorname{diag}}
\def\prop{\operatorname{prop}}
\def\diam{\operatorname{diam}}
\def\supp{\operatorname{supp}}

\author[L.~Guo]{Liang Guo}
\address[L.~Guo]{Centre for Mathematical Sciences, Lund University, Box 118, Lund, 22100, Sweden}
\email{liang.guo@math.lth.se}

\author[K.~Li]{Kang Li}
\address[K.~Li]{Centre for Mathematical Sciences, Lund University, Box 118, Lund, 22100, Sweden}
\email{kang.li@math.lth.se}

\author[J.~Zhang]{Jiawen Zhang}
\address[J.~Zhang]{School of Mathematical Sciences, Fudan University, 220 Handan Road,
Shanghai, 200433, China.}
\email{jiawenzhang@fudan.edu.cn}

\author[J.~Zhu]{Jingming Zhu}
\address[J.~Zhu]{College of Data Science, Jiaxing University, 899 Guangqiong Road, Jiaxing, 314000, China.}
\email{jingmingzhu@zjxu.edu.cn}

\title[Nuclear dimension of uniform Roe algebras]{Nuclear dimension of uniform Roe algebras}
\date{\today}
\begin{document}

\begin{abstract}
We prove that for every discrete metric space $X$ of bounded geometry with Yu's Property~A, the uniform Roe algebra $C^*_u(X)$ has nuclear dimension at most one. This answers in the negative a long-standing question of Winter and Zacharias on whether the nuclear dimension of $C^*_u(X)$ always equals the asymptotic dimension of $X$. Consequently, $\dimnuc C^*_u(X)\in\{0,1,\infty\}$ for every discrete metric space $X$, with or without bounded geometry.
\end{abstract}

\maketitle
\tableofcontents

\section{Introduction}\label{sec: intro}

The nuclear dimension of $C^*$-algebras was introduced by Winter and Zacharias \cite{WZ10} as a noncommutative analogue of the covering dimension of topological spaces. It refines the decomposition rank of Kirchberg and Winter \cite{KW04}, and recovers the covering dimension of $Y$ for a commutative $C^*$-algebra $C_0(Y)$. Since then, finite nuclear dimension has become one of the central regularity properties in the classification program of simple nuclear $C^*$-algebras. We refer to \cite{STW} for a survey. A striking feature is that the nuclear dimension of simple $C^*$-algebras can only take the values $0$, $1$ and $\infty$ \cite{Win12, CE20, CETWW21}.

Uniform Roe algebras form another natural class of examples coming from coarse geometry, and they are not simple unless the underlying space is finite. For a discrete metric space $X$ of bounded geometry, the uniform Roe algebra $C^*_u(X)$ is the norm closure of the bounded operators on $\ell^2(X)$ of finite propagation. Winter and Zacharias showed that $\dimnuc C^*_u(X)\le\asdim X$ \cite[Theorem~8.5]{WZ10}, where $\asdim X$ denotes the asymptotic dimension of $X$. This suggests that the nuclear dimension of $C^*_u(X)$ might recover the asymptotic dimension of $X$, as in the commutative case, and Winter and Zacharias asked whether $\dimnuc C^*_u(X)=\asdim X$ \cite[Question~9.5]{WZ10}. This question was recently listed again in \cite[Problem~LXXXVII]{STW}, where even $X=\IZ^2$ is recorded as open.

Several results seem to support a positive answer. Uniform Roe algebras are rigid: if two spaces of bounded geometry have isomorphic uniform Roe algebras, then they are coarsely equivalent. This was first proved by \v{S}pakula and Willett in \cite{SW13} under Yu's Property~A \cite{Yu00}, and by Baudier, Braga, Farah, Khukhro, Vignati and Willett in general \cite{BBFKVW22}. Moreover, Vignati showed that the two spaces are even bijectively coarsely equivalent \cite{Vignati2026}. Since the asymptotic dimension is a coarse invariant, it follows that $\asdim X$ is determined by the $C^*$-algebra $C^*_u(X)$, and one may expect it to be read off from some invariant of $C^*_u(X)$. This is indeed the case for dimension zero: Li and Willett showed that $\dimnuc C^*_u(X)=0$ if and only if $\asdim X=0$ \cite[Corollary~1.5]{LW18}.

More recently, Li, Liao and Winter introduced the notion of diagonal dimension, a version of the nuclear dimension which also takes a distinguished commutative subalgebra into account \cite{LLW}. They showed that the diagonal dimension of $\ell^\infty(X)\subseteq C^*_u(X)$ is exactly $\asdim X$. 
All these results suggest that $\dimnuc C^*_u(X)$ and $\asdim X$ should coincide.

However, in this paper we show that the answer to the question of Winter and Zacharias is negative. Our main result is the following.

\begin{Thm}\label{thm: main theorem}
Let $X$ be a countable discrete metric space of bounded geometry. If $X$ has Property~A, then $\dimnuc C^*_u(X)\le1$.
\end{Thm}

In particular, $\dimnuc C^*_u(\IZ^d)=1$ for every $d\ge1$, while $\asdim\IZ^d=d$. Theorem~\ref{thm: main theorem} also applies to spaces of infinite asymptotic dimension with Property~A, such as the lamplighter group $\IZ\wr\IZ$. On the other hand, $C^*_u(X)$ is nuclear if and only if $X$ has Property~A \cite{STY02, BO08}, so $\dimnuc C^*_u(X)=\infty$ whenever Property~A fails. Combining these results, we obtain the following trichotomy. 

\begin{Cor}\label{cor: trichotomy}
Let $X$ be a countable discrete metric space of bounded geometry. Then
\[
  \dimnuc C^*_u(X)=
  \begin{cases}
    0, & \asdim X=0,\\
    1, & X\text{ has Property A and }\asdim X>0,\\
    \infty, & X\text{ does not have Property A}.
  \end{cases}
\]
\end{Cor}

Let us briefly explain the idea of the proof. To prove $\dimnuc C^*_u(X)\le\asdim X$, Winter and Zacharias cut operators by a slowly varying partition of unity on $X$ whose functions are supported in the members of $\asdim X+1$ uniformly bounded families of subsets. With only two colors, such a partition of unity on $X$ does not exist in general. Instead, we use a noncommutative version, which we call a \emph{twisted partition of unity}. More precisely, we write $1=h+V^*V$, where $h$ is a function on $X$ and $V\in C^*_u(X)$ is an operator such that $VV^*$ is close to a function $g$ on $X$, and both $h$ and $g$ cut operators into block diagonal pieces. In Section~\ref{sec: twisted partition of unity}, we show that such twisted partitions of unity always exist for spaces with Property~A. The construction of $V$ relies on a purely coarse geometric construction, which we call \emph{geometric matching}, and which we establish in Section~\ref{sec: geometric matching}. Finally, in Section~\ref{sec: proof of the main theorem}, we deduce Theorem~\ref{thm: main theorem} from the twisted partition of unity, and prove Corollary~\ref{cor: trichotomy}.

\subsection*{Acknowledgement}

Jiawen Zhang was partly supported by the National Key R{\&}D Program of China 2022YFA100700 and NSFC 12422107.

\subsection*{AI Declaration}

A first proof of Theorem~\ref{thm: main theorem} was obtained with the assistance of GPT-6.1 Sol, which relied on the Lov\'asz Local Lemma. By working out the case $X=\IZ^2$ by hand, the authors observed that the randomness is only used to find far apart points of the same local type near the boundaries of a partition, and that such points can be found deterministically. The deterministic proof presented here was developed from this observation. AI tools were used to work out and simplify the argument and for language polishing. The authors checked every argument and determined the structure and presentation of the paper, and take full responsibility for its content.

\section{Geometric matching}\label{sec: geometric matching}

Let $(X, d)$ be a countable discrete metric space of bounded geometry. For $R > 0$, we denote by $B_R(x)$ the closed ball of radius $R$ centered at $x \in X$. For any $F\subseteq X$, we write $B_R(F)=\bigcup_{x\in F} B_R(x)$. By \emph{bounded geometry}, we mean that there exists a function $N: \mathbb{R}^+ \to \mathbb{N}$ such that for all $x \in X$ and $R > 0$, the ball $B_R(x)$ contains at most $N(R)$ elements. Throughout this paper, when we say \emph{$X$ is a space}, we mean $X$ is a metric space with bounded geometry.

In this paper, we shall also consider colored metric spaces. A \emph{finite coloring} on $X$ is a map $c: X\to\{1,\cdots,\kappa\}$. A \emph{partition} $\P$ of $X$ is a family of nonempty, uniformly bounded, pairwise disjoint subsets of $X$ such that their union is $X$. A set in $\P$ is called a \emph{block}. For any $x\in X$, we shall denote by $\P(x)$ the unique block containing $x$. A partition $\Q$ is coarser than $\P$ if every block in $\Q$ can be writen as a union of blocks of $\P$. For any $r\geq 0$, we define the $r$-inner boundary of $\P$ to be
$$\partial_r\P=\{x\in X \mid B_r(x)\nsubseteq \P(x)\}.$$
One checks directly that if $\Q$ is coarser than $\P$, then $\P(x)\subseteq\Q(x)$ and hence $\partial_r\Q\subseteq \partial_r\P$.

Our goal for this section is to prove the following geometric matching lemma.

\begin{Pro}\label{pro: geometric matching}
Let $X$ be a space with Property A with a finite color $c: X\to\{1,\cdots,\kappa\}$. For any $m\in\IN$ and $r>0$, there exists a partition $\P$ of $X$ and an injection
$$\mu:\partial_r\P\times\{1, \cdots, m\} \to X$$
such that for any $x,x'\in \partial_r\P$ and $j,j'\in\{1,\cdots, m\}$:\begin{enumerate}[(1)]
\item $c(\mu(x,j)) = c(x)$ and $\sup_{x,j}d(x,\mu(x,j))  < \infty$;
\item $B_r(\mu(x,j)) \subseteq \P(\mu(x,j))$;
\item $d(\mu(x,j), \mu(x',j')) > 2r$ whenever $(x,j) \ne (x',j')$. 
\end{enumerate}
\end{Pro}

In words, the proposition says that for any given $m$ and $r$ one can find a partition $\P$ such that every point $x\in\partial_r\P$ is assigned $m$ \emph{partners} $\mu(x,1),\dots,\mu(x,m)$ of the same color as $x$. These partners lie within a uniformly bounded distance of $x$, sit deep inside the blocks of $\P$, and are well separated from one another.

We begin with a brief review of Property~A. For further details, we refer the reader to \cite{Wil09}. 

\begin{Def}\label{def: property A}
A space $(X, d)$ has \emph{Property~A} if for every $R > 0$ and $\varepsilon > 0$, there exist $S > 0$ and a map $X \to \ell^{1}(X)$, $x \mapsto \xi_x$, such that
\begin{enumerate}[(i)]
\item $\|\xi_x\|_1 = 1$ and $\xi_x(y) \geq 0$ for all $x, y \in X$;
\item $\supp(\xi_x) \subseteq B_S(x)$ for all $x \in X$;
\item $\|\xi_x - \xi_y\|_1 < \varepsilon$ whenever $d(x, y) < R$.
\end{enumerate}
\end{Def}

There are many equivalent characterizations of Property~A. In this paper, we shall also use the \emph{metric sparsification property}, introduced by X.~Chen, R.~Tessera, X.~Wang, and G.~Yu in \cite{CTWY08}. For convenience, we call a finite family $\U = \{U_1, \dots, U_n\}$ of nonempty subsets of $X$ an \emph{$(R,S)$-packing} if
\begin{enumerate}[(i)]
\item $\diam(U_i) \leq S$ for every $i = 1, \dots, n$;
\item $d(U_i, U_j) > R$ for all $i \neq j$.
\end{enumerate}
For any packing $\U = \{U_1, \dots, U_n\}$, we shall write $\bigsqcup\U =  U_1 \sqcup \cdots \sqcup U_n$.

\begin{Def}\label{def: MSP}
A space $X$ has the \emph{metric sparsification property (MSP)} if for every $\varepsilon > 0$ and $R > 0$, there exists $S > 0$ such that for every finite positive measure $\mu$ on $X$, there exists an $(R,S)$-packing $\U$ satisfying
\[
\mu(\textstyle\bigsqcup\U) \geq (1 - \varepsilon)\mu(X).
\]
\end{Def}

For a countable discrete space $X$, every subset is measurable, and a finite positive measure $\mu$ can be identified with a function $\mu: X \to [0, \infty)$ satisfying $\sum_{x \in X} \mu(x) < \infty$. Notice that the definition of the MSP used here differs slightly from the original one in \cite[Definition~3.1]{CTWY08}: instead of fixing a constant $c \in (0, 1)$, we require the inequality to hold for any $c = 1 - \varepsilon$ arbitrarily close to $1$. The equivalence of these two formulations is shown in \cite[Proposition~3.3]{CTWY08}. Another minor difference is that we restrict ourselves to finite families $\U$. Since $\mu$ is finite, one can always truncate an infinite family by discarding a tail of sufficiently small measure to obtain a finite subfamily satisfying the desired inequality. It is proved by \cite{BNSWW13} and \cite{Sak14} that MSP is equivalent to Property A.

\begin{Lem}\label{lem: volume expanding}
Let $X$ be a space with Property A. For any $R>0$ and $k\in\IN$, there exist $S>0$ and $D_0>0$ such that for any finite subset $Y\subseteq X$, there exists an $(R,S)$-packing $\U$ in $Y$ such that for any $F\subseteq Y\setminus \bigsqcup \U$, 
$$\#(Y\cap B_{D_0}(F)) \geq k\cdot \# F.$$
\end{Lem}

\begin{proof}
Let $S$ be given by Definition~\ref{def: MSP} for $R$ and $\varepsilon=1/(4k)$, and put $D_0=4k(R+S)$.

Let $Y\subseteq X$ be finite. Every member of an $(R,S)$-packing in $Y$ is a subset of $Y$, so there are only finitely many such packings. Choose one $\U$ such that $\#\bigsqcup\U$ is the largest.

Assume for a contradiction that there exists $F\subseteq Y\setminus \bigsqcup\U$ satisfies $\#(Y\cap B_{D_0}(F))<k\cdot\#F$. Then $F$ is non-empty. For $0\le j\le 4k$, put
$$F_j=Y\cap B_{j(R+S)}(F).$$
Then $F=F_0\subseteq F_1\subseteq\cdots\subseteq F_{4k}=Y\cap B_{D_0}(F)$. By assumption $\#F_{4k}<k\cdot\#F$. The sets $F_{j+1}\setminus F_j$ are pairwise disjoint subsets of $F_{4k}$, hence $\sum_{j=0}^{4k-1}\#(F_{j+1}\setminus F_{j})<k\cdot\#F$. Then there exists some $j$
satisfies
$$\#(F_{j+1}\setminus F_{j})<\tfrac{\#F}4 .$$

By metric sparsification, considering the counting measure on $F_j$, one can find an $(R,S)$-packing $\mathcal V$ in $F_j$ with $\#\bigsqcup\V \ge (1-\frac1{4k})\,\#F_j$. Define
$$\U'=\V\cup \{U\in \U \mid d(U,\textstyle\bigsqcup\V) > R \}\quad\text{and}\quad\U_{\rm deleted}=\{U\in \U \mid d(U,\textstyle\bigsqcup\V) \leq R \}$$
Thus, $\U'$ is again an $(R,S)$-packing in $Y$. Each $U\in \U_{\rm deleted}$ has diameter at most $S$ and contains a point within distance $R$ of $V\subseteq F_j\subseteq B_{j(R+S)}(F)$ for some $V\in\V$. Hence $U\subseteq F_{j+1}$. Therefore $\bigsqcup\U_{\rm deleted}\subseteq F_{j+1}\setminus F$,
and
$$\#\textstyle\bigsqcup\U_{\rm deleted}\le\#F_{j+1}-\#F = \#F_j-\#F+\#(F_{j+1}\setminus F_j).$$
By assumption, $\#F_j\leq \#F_{4k} < k\cdot\#F$, we then conclude that
$$\#\textstyle\bigsqcup\U'-\#\textstyle\bigsqcup\U = \#\textstyle\bigsqcup\V-\#\textstyle\bigsqcup\U_{\rm deleted} \ge \#F-\frac{\#F_j}{4k}-\#(F_{j+1}\setminus F_j) > \#F-\frac{\#F}{4}-\frac{\#F}{4} > 0,$$
which contradicts the maximality of $\#\textstyle\bigsqcup\U$.
\end{proof}

\begin{Lem}\label{lem: bounded step}
Let $X$ be a space with Property A. For any $R>0$, $b'\ge 0$ and $k\in\IN$, there exist $b>0$ and $D_0>0$ such that for any partition $\Q'$ of $X$ with $\sup_{D\in\Q'}\diam(D)\le b'$ and any finite subset $Y\subseteq X$, there exist a partition $\Q$ coarser than $\Q'$ with $\sup_{D\in\Q}\diam(D)\le b$ and an injective map
\[\nu\colon (Y\cap\partial_R\Q)\times\{1,\dots,k\}\to Y\]
such that $d(x,\nu(x,j))\le D_0$ for all $x\in Y\cap \partial_R\Q$ and $j\in\{1,\dots,k\}$.
\end{Lem}

To prove this lemma, we first recall Hall's theorem.

\begin{Lem}[\cite{Hal35}]\label{lem: Hall}
Let $B$ be a set, let $I$ be a finite set, and for each $i\in I$ let $\Gamma(i)\subseteq B$ be a finite subset. If for every $I'\subseteq I$,
\[\#\Bigl(\bigcup_{i\in I'}\Gamma(i)\Bigr)\ \ge\ \#I'\]
then there is an injective map $\phi\colon I\to B$ such that
$\phi(i)\in\Gamma(i)$ for every $i\in I$.\qed
\end{Lem}

\begin{proof}[Proof of Lemma~\ref{lem: bounded step}]
Let $S$ and $D_0$ be given by Lemma~\ref{lem: volume expanding} for $2R+b'$ and $k$, and put
$b=S+2R+2b'$. Let $\Q'$ and $Y$ be given, and let $\U$ be the $(2R+b',S)$-packing in $Y$
provided by Lemma~\ref{lem: volume expanding}. For $U\in\U$, let
\[
  Q_U=\bigcup\{D\in\Q' \mid D\cap B_R(U)\neq\emptyset\}.
\]
If a block $D\in\Q'$ met both $B_R(U)$ and $B_R(U')$ for distinct $U,U'\in\U$, then picking $u\in D\cap B_R(U)$ and $u'\in D\cap B_R(U')$ we would get
\[
  d(U,U')\le R+d(u,u')+R\le 2R+b',
\]
contradicting $d(U,U')>2R+b'$. Hence the family of sets $\{Q_U\}_{U\in\U}$ are pairwise disjoint. Every point of $Q_U$ lies in a block of diameter at most $b'$ that meets $B_R(U)$, hence is within distance $b'+R$ of $U$. Since $\diam(U)\le S$, we get $\diam(Q_U)\le (b'+R)+S+(R+b')=b$. Let $\Q$ consist of the sets $\{Q_U\}_{U\in\U}$ together with the blocks of $\Q'$ that meet no $B_R(U)$, which have diameter at most $b'\le b$.  Thus, $\Q$ is a partition of $X$ coarser than $\Q'$.

Let $y\in U$ for some $U\in\U$. Each point $z\in B_R(y)\subseteq B_R(U)$ lies in a block of $\Q'$ that meets $B_R(U)$, so $B_R(y)\subseteq Q_U=\Q(y)$. Thus $y\notin\partial_R\Q$. As a result, $Y\cap\partial_R\Q\subseteq Y\setminus\bigsqcup\U$.

Finally, we apply Lemma~\ref{lem: Hall} with $B=Y$, $I=(Y\cap\partial_R\Q)\times\{1,\dots,k\}$ and $\Gamma(x,j)=Y\cap B_{D_0}(x)$. Let $I'\subseteq I$, and let $F\subseteq Y$ be the set of first coordinates of elements of $I'$. Then $F\subseteq Y\setminus\bigsqcup\U$ and $\#I'\le k\cdot\#F$. By Lemma~\ref{lem: volume expanding}
\[
  \#\Bigl(\bigcup_{i\in I'}\Gamma(i)\Bigr)=\#\bigl(Y\cap B_{D_0}(F)\bigr)
  \ge k\cdot\#F\ge\#I'.
\]
Hall's theorem yields an injective map $\nu\colon I\to Y$ with $\nu(x,j)\in Y\cap B_{D_0}(x)$, that is, $d(x,\nu(x,j))\le D_0$.
\end{proof}

We now prove a \emph{finite set version} of Proposition~\ref{pro: geometric matching} and show that the constants do not depend on the finite set.

\begin{Pro}\label{pro: finite matching}
Let $X$ be a space with Property A and let $c\colon X\to\{1,\dots,\kappa\}$ be a finite coloring. For any $m\in\IN$ and $r>0$, there exist $b>0$ and $D>0$ such that for any \emph{finite} subset $W\subseteq X$, there exist a partition $\P$ of $X$ whose blocks have diameter at most $b$ and an injective map
\[
  \mu\colon (W\cap\partial_r\P)\times\{1,\dots,m\}\to X
\]
such that for any $x,x'\in W\cap\partial_r\P$ and $j,j'\in\{1,\dots,m\}$:
\begin{enumerate}[(1)]
\item $c(\mu(x,j))=c(x)$ and $d(x,\mu(x,j))\le D$;
\item $B_r(\mu(x,j))\subseteq\P(\mu(x,j))$;
\item $d(\mu(x,j),\mu(x',j'))>2r$ whenever $(x,j)\neq(x',j')$.
\end{enumerate}
\end{Pro}

The point is that $b$ and $D$ do not depend on $W$. To prove this proposition, we need the following counting lemma.

\begin{Lem}\label{lem: separated choice}
Let $X$ be a space. For any $R\ge0$, there exists $K\in\IN$ such that for any pairwise disjoint finite subsets $L_1,\dots,L_n\subseteq X$ with $\#L_i\ge K$ for every $i$, there exist $y_i\in L_i$ such that $d(y_i,y_l)>R$ for all $i\neq l$.
\end{Lem}

\begin{Rem}
Lemma~\ref{lem: separated choice} is a special case of the following theorem of Haxell \cite{Hax01}: if $G$ is a finite graph of maximum degree at most $\Delta$ whose vertex set is partitioned into sets $V_1,\dots,V_n$ of size at least $2\Delta$, then $G$ has an independent set containing exactly one vertex of each $V_i$. Indeed, applying it to the graph on $L_1\sqcup\cdots\sqcup L_n$ in which two distinct points are adjacent if their distance is at most $R$, whose maximum degree is at most $N(R)-1$, one may even take $K=2(N(R)-1)$. The reader willing to accept Haxell's theorem may skip the following proof, which is an elementary counting argument independent of \cite{Hax01}, included to keep the paper self-contained.
\end{Rem}

\begin{proof}
Put $K=4N(R)$. Removing points if necessary, we may assume $\#L_i=K$ for every $i$. We may also assume $n\ge2$, the case $n=1$ being trivial. For a nonempty subset $I\subseteq\{1,\dots,n\}$, we call a family $(y_i)_{i\in I}$ with $y_i\in L_i$ \emph{admissible} if $d(y_i,y_l)>R$ for all distinct $i,l\in I$, and we denote by $\alpha(I)$ the number of admissible families. Note that $\alpha(\{i\})=\#L_i=K$ for every $i$. We show by induction on $\#I\ge1$ that
\[
  \alpha(I\cup\{l\})\ \ge\ \frac K2\cdot\alpha(I)\qquad\text{for every } l\notin I.
\]

Suppose first that $\#I=1$, say $I=\{i\}$. There are $K^2$ pairs $(y_i,y)$ with $y_i\in L_i$ and $y\in L_l$. For each $y\in L_l$, at most $\#B_R(y)\le N(R)=K/4$ points $y_i\in L_i$ satisfy $d(y_i,y)\le R$. Hence at most $K^2/4$ pairs are not admissible, and
\[
  \alpha(\{i,l\})\ \ge\ K^2-\frac{K^2}4=\frac{3K^2}4\ \ge\ \frac K2\cdot K
  =\frac K2\cdot\alpha(\{i\}).
\]

Let $\#I\ge2$, and assume the claim holds for all nonempty subsets with fewer than $\#I$ elements. A pair consisting of an admissible family $(y_i)_{i\in I}$ and a point $y\in L_l$ is an admissible family on $I\cup\{l\}$ unless $d(y,y_i)\le R$ for some $i\in I$. Since the sets $L_i$ are pairwise disjoint, for each $y\in L_l$, there are at most $\#B_R(y)\le K/4$ pairs $(i,z)$ with $i\in I$, $z\in L_i$ and $d(z,y)\le R$. Fix such a pair $(i,z)$. An admissible family on $I$ with $y_i=z$ is determined by its restriction to $I\setminus\{i\}$, which is again admissible. Since $I\setminus \{i\}$ is nonempty, the induction hypothesis gives
\[
  \alpha(I\setminus\{i\})\ \le\ \frac2K\cdot\alpha(I).
\]
So for each $y\in L_l$, at most $\frac K4\cdot\frac2K\cdot\alpha(I)=\alpha(I)/2$ admissible families on $I$ fail to combine with $y$. Hence at most $K\cdot\alpha(I)/2$ of the $K\cdot\alpha(I)$ pairs are not admissible, and
\[
  \alpha(I\cup\{l\})\ \ge\ \frac K2\cdot\alpha(I).
\]

This proves the claim. Applying it repeatedly, starting from $\{1\}$, we obtain $\alpha(\{1,\dots,n\})\ge(K/2)^{n-1}K>0$, so an admissible family on $\{1,\dots,n\}$ exists.
\end{proof}

\begin{proof}[Proof of Proposition~\ref{pro: finite matching}]
Let $W\subseteq X$ be finite. We shall prove this proposition in the following steps.

\emph{Step 1: Find separated partners.} Let $K$ be given by Lemma~\ref{lem: separated choice} for $2r$, and let
$\Q_0=\bigl\{\{x\}\mid x\in X\bigr\}$ be the partition of $X$ into one-point sets. We construct partitions $\Q_0,\Q_1,\dots,\Q_\kappa$ of $X$ inductively, each coarser than the previous one, the $i$-th step takes care of the points of color $i$. Suppose that $\Q_{i-1}$ has been constructed for some $1\le i\le\kappa$. Apply Lemma~\ref{lem: bounded step} to $2r$, $mK$, $\Q'=\Q_{i-1}$ and the finite set $Y=W\cap c^{-1}(i)$, we thus obtain a partition $\Q_i$ of $X$ coarser than $\Q_{i-1}$ together with an injective map
\[
  \nu_i\colon \bigl(W\cap c^{-1}(i)\cap\partial_{2r}\Q_i\bigr)\times\{1,\dots,mK\} \to W\cap c^{-1}(i).
\]
We emphasize that each $\Q_i$ is a partition of the whole space $X$, the set $W\cap c^{-1}(i)$ only determines which boundary points receive candidates in the $i$-th step. The constants in Lemma~\ref{lem: bounded step} depend only on $X$, $2r$, $mK$ and the diameter bound of the blocks of $\Q_{i-1}$. Hence, by induction on $i$, there exist $b'>0$ and $D>0$ depending only on $X$, $\kappa$, $m$ and $r$ such that the blocks of every $\Q_i$ have diameter at most $b'$ and every $\nu_i$ moves points by at most $D$.

Since $\Q_{i+1},\dots,\Q_\kappa$ are obtained from $\Q_i$ by further coarsening, they create no new boundary points as $\partial_{2r}\Q_\kappa\subseteq\partial_{2r}\Q_i$. Hence for every $x\in W\cap\partial_{2r}\Q_\kappa$ and $1\le l\le mK$, the point $\nu_{c(x)}(x,l)$ is defined, has color $c(x)$, and lies within distance $D$ of $x$. Moreover, these points are distinct for distinct $(x,l)$. Indeed, if $\nu_{c(x)}(x,l)=\nu_{c(x')}(x',l')$, then $c(x)=c(x')$ and the injectivity of $\nu_{c(x)}$ gives $(x,l)=(x',l')$.

For $x\in W\cap\partial_{2r}\Q_\kappa$ and $j\in\{1,\dots,m\}$, put
\[
  L_{x,j}=\{\nu_{c(x)}(x,l)\mid (j-1)K<l\le jK\}.
\]
These sets are pairwise disjoint and each has $K$ elements. By Lemma~\ref{lem: separated choice}, there exist $y_{x,j}\in L_{x,j}$ such that $d(y_{x,j},y_{x',j'})>2r$ whenever $(x,j)\neq(x',j')$.

\emph{Step 2: Rerouting the boundary.} By Step 1, the balls $B_r(y_{x,j})$ are pairwise disjoint, so we may define a map $\pi\colon X\to X$ by
\[
  \pi(z)=\begin{cases}
    y_{x,j}, & \text{if } z\in B_r(y_{x,j}) \text{ for some } (x,j),\\
    z, & \text{otherwise.}
  \end{cases}
\]
Then $d(z,\pi(z))\le r$ for all $z\in X$. Let $\P$ be the partition of $X$ consisting of the nonempty sets $\pi^{-1}(E)$, $E\in\Q_\kappa$. In other words, $\P(z)=\pi^{-1}\bigl(\Q_\kappa(\pi(z))\bigr)$.

If $z'\in\P(z)$, then $\pi(z)$ and $\pi(z')$ lie in the same block of $\Q_\kappa$, so
\[
  d(z,z')\le d(z,\pi(z))+d(\pi(z),\pi(z'))+d(\pi(z'),z')\le b'+2r .
\]
Since $\pi$ is constant on $B_r(y_{x,j})$ with value $y_{x,j}=\pi(y_{x,j})$, we have $B_r(y_{x,j})\subseteq\P(y_{x,j})$. We claim that
\[
  \partial_r\P\subseteq\partial_{2r}\Q_\kappa .
\]
Indeed, let $x\notin\partial_{2r}\Q_\kappa$, that is, $B_{2r}(x)\subseteq\Q_\kappa(x)$. For any $z\in B_r(x)$, $d(x,\pi(x))\le r$ and $d(x,\pi(z))\le d(x,z)+d(z,\pi(z))\le 2r$, so $\pi(x),\pi(z)\in B_{2r}(x)\subseteq\Q_\kappa(x)$. Hence $z\in\P(x)$. As a result, $B_r(x)\subseteq\P(x)$ and $x\notin\partial_r\P$.

\emph{Conclusion.} By the claim, every $x\in W\cap\partial_r\P$ lies in $W\cap\partial_{2r}\Q_\kappa$, so we may define $\mu(x,j)=y_{x,j}$. Condition (1) and (3) holds by Step 1, condition (2) by Step 2. In particular, $\mu$ is injective. Finally, the blocks of $\P$ have diameter at most $b=b'+2r$, and both $b$ and $D$ depend only on $X$, $\kappa$, $m$ and $r$, not on $W$.
\end{proof}

Now, we are ready to prove the geometric matching for the general case.

\begin{proof}[Proof of Proposition~\ref{pro: geometric matching}]
Let $b$ and $D$ be given by Proposition~\ref{pro: finite matching}. Since $X$ is countable, we may enumerate $X=\{x_1,x_2,\dots\}$. For each $n\in\IN$, we apply Proposition~\ref{pro: finite matching} to the finite set $W=\{x_1,\dots,x_n\}$ and obtain a partition $\P_n$ of $X$ whose blocks have diameter at most $b$, together with an injective map
\[
  \mu_n\colon \bigl(\{x_1,\dots,x_n\}\cap\partial_r\P_n\bigr)\times\{1,\dots,m\}\to X
\]
satisfying conditions (1)--(3) of Proposition~\ref{pro: finite matching}.

\emph{The limit partition.} For every $x\in X$, the block $\P_n(x)$ has diameter at most $b$, so it is a subset of the finite ball $B_b(x)$. Hence $\P_n(x)$ takes only finitely many values as $n$ varies. We use a diagonal argument to pass to a subsequence along which all these blocks stabilize. More precisely, since $\P_n(x_1)$ takes only finitely many values, there is an infinite set $N_1\subseteq\IN$ on which $\P_n(x_1)$ is constant. Inductively, there are infinite sets $N_1\supseteq N_2\supseteq\cdots$ such that $\P_n(x_1),\dots,\P_n(x_k)$ are constant for $n\in N_k$. Choose $n_1<n_2<\cdots$ with $n_k\in N_k$, and define
\[
  \P(x_i)=\P_{n_i}(x_i)\qquad\text{for every } i .
\]
Since $n_k\in N_k\subseteq N_i$ for $k\ge i$, we have $\P_{n_k}(x_i)=\P(x_i)$ for all $k\ge i$.

We stress that the nestedness $N_1\supseteq N_2\supseteq\cdots$ is essential. For any $i,l$ and any $k\ge\max\{i,l\}$, we have $n_k\in N_i\cap N_l$, so $\P(x_i)$ and $\P(x_l)$ are blocks of one and the same partition $\P_{n_k}$. Had we chosen the blocks of different points independently, it could happen that $x_l$ lies in the block chosen for $x_i$ but is assigned a different block.

The family $\{\P(x)\}_{x\in X}$ forms a partition of $X$. Indeed, $x_i\in\P_{n_i}(x_i)=\P(x_i)$ for every $i$. Now let $y\in\P(x)$, say $x=x_i$ and $y=x_l$, and choose $k\ge\max\{i,l\}$. Then $y\in\P(x)=\P_{n_k}(x)$, so $\P_{n_k}(y)=\P_{n_k}(x)$ because $\P_{n_k}$ is a partition. Since also $\P_{n_k}(y)=\P(y)$ and $\P_{n_k}(x)=\P(x)$, we get $\P(y)=\P(x)$. Consequently, if $z\in\P(x)\cap\P(y)$, then $\P(x)=\P(z)=\P(y)$; that is, two sets $\P(x)$ and $\P(y)$ are either equal or disjoint. Finally, each $\P(x_i)=\P_{n_i}(x_i)$ is a block of $\P_{n_i}$, so it has diameter at most $b$.

\emph{The limit partners.} Let $x\in\partial_r\P$, say $x=x_i$, and let $j\in\{1,\dots,m\}$. For every $k\ge i$, we have $n_k\ge k\ge i$, so $x\in\{x_1,\dots,x_{n_k}\}$. Moreover, $\P_{n_k}(x)=\P(x)$, so $x\in\partial_r\P_{n_k}$. Hence $\mu_{n_k}(x,j)$ is defined for all $k\ge i$, and it lies in the finite ball $B_D(x)$. Since $\partial_r\P\times\{1,\dots,m\}$ is countable, the same diagonal argument as above allows us to pass to a subsequence $(n_{k_l})_l$ of $(n_k)_k$ such that, for every $x\in\partial_r\P$ and $j$, the point $\mu_{n_{k_l}}(x,j)$ is constant for all sufficiently large $l$. We denote this eventual value by $\mu(x,j)$.

Finally, we verify that $\P$ and $\mu$ satisfy conditions (1)--(3). Let $x,x'\in\partial_r\P$ and $j,j'\in\{1,\dots,m\}$, and choose $l$ so large that, with $n=n_{k_l}$, we have $\mu_n(x,j)=\mu(x,j)$, $\mu_n(x',j')=\mu(x',j')$ and
$\P_n(\mu(x,j))=\P(\mu(x,j))$.
\begin{enumerate}[(1)]
\item By condition (1) for $\mu_n$, we have $c(\mu(x,j))=c(x)$ and $d(x,\mu(x,j))\le D$. In particular, $\sup_{x,j}d(x,\mu(x,j))\le D<\infty$.
\item By condition (2) for $\mu_n$, $B_r(\mu(x,j))\subseteq\P_n(\mu(x,j))=\P(\mu(x,j))$.
\item If $(x,j)\neq(x',j')$, then condition (3) for $\mu_n$ gives $d(\mu(x,j),\mu(x',j'))>2r$.
\end{enumerate}
In particular, $\mu$ is injective. This completes the proof.
\end{proof}

\section{Twisted partition of unity}\label{sec: twisted partition of unity}

In this section, we shall introduce a notion of \emph{twisted partition of unity} and prove that every space with Property A admits a twisted partition of unity.

Let $\{\delta_x\}_{x\in X}$ be the standard basis of $\ell^2(X)$. For any $T\in\B(\ell^2(X))$ and $x,y\in X$, we write $T_{xy}=\langle T\delta_y,\delta_x\rangle$, and the \emph{propagation} of $T$ is defined by
$$\prop(T)=\sup\{d(x,y)\mid T_{xy}\neq0\}.$$
The \emph{uniform Roe algebra} $C^*_u(X)$ is the norm closure of the set of bounded operators on $\ell^2(X)$ of finite propagation.
A bounded function $f\colon X\to\IC$
is identified with the multiplication operator $\delta_x\mapsto f(x)\delta_x$.
For $E\subseteq X$, $\chi_E$ denotes the multiplication operator by the
characteristic function of $E$, i.e., the orthogonal projection onto
$\ell^2(E)$.
Given a partition $\P$ of $X$, an operator $T\in \B(\ell^2(X))$ is called \emph{$\P$-diagonal} if for any distinct $P,P'\in\P$, one has that $\chi_{P}T\chi_{P'}=0$. 

\begin{Def}\label{def: twisted partition of unity}
Let $\P$ be a partition of $X$. For any finite set $\F\subseteq C^*_u(X)$ and $\varepsilon>0$, a \emph{twisted partition of unity subordinate to $\P$} (for $\F$, $\varepsilon$) consists of two functions $h,g\colon X\to[0,1]$ and an operator $V\in C^*_u(X)$ such that
\begin{enumerate}[(1)]
\item $V^*V=1-h$;
\item $\|VV^*-g\|\le\varepsilon$;
\item $h^{1/2}Th^{1/2}$ and $g^{1/2}Tg^{1/2}$ are $\P$-diagonal for every $T\in\F$;
\item $\|[h^{1/2},T]\|<\varepsilon$ and $\|[V,T]\|<\varepsilon$ for every $T\in\F$.
\end{enumerate}
\end{Def}

By (1), we have that $h+V^*V=1$ is indeed a partition of unity. The operator $(1-h)^{1/2}T(1-h)^{1/2}$ is in general not $\P$-diagonal. Condition (2) says that the ``twist'' $VV^*$ of $1-h=V^*V$ is close to the function $g$, and by
(3) the compression by $g^{1/2}$ is again $\P$-diagonal.

\begin{Pro}\label{pro: twisted partition of unity}
Let $X$ be a space with Property A. Let $\F\subseteq C^*_u(X)$ be a finite set of contractions with $\F^*=\F$ and $\prop(T) < \infty$ for every $T\in\F$. For any $\varepsilon > 0$, there exist a partition $\P$ of $X$ and a twisted partition of unity subordinate to $\P$ for $(\F,\varepsilon)$.
\end{Pro}

We shall use the rest of this section to prove Proposition~\ref{pro: twisted partition of unity}. Before we get into the detailed construction, we need the following two lemmas.

\begin{Lem}\label{lem: local types}
Let $X$ be a space and let $\F\subseteq\B(\ell^2(X))$ be a finite set of contractions. For any $0\le S\le L$ and $\varepsilon>0$, there exists a finite coloring $c\colon X\to\{1,\dots,\kappa\}$ such that, whenever $c(z)=c(y)$, there is a bijection $\varphi\colon B_L(z)\to B_L(y)$ mapping $B_S(z)$ onto $B_S(y)$ whose induced partial isometry $v$, given by $v\delta_x=\delta_{\varphi(x)}$ for $x\in B_L(z)$ and $v\delta_x=0$ otherwise, satisfies that for any $T\in \F$,
\[
  \bigl\|\chi_{B_L(z)}T\chi_{B_L(z)}-v^*\chi_{B_L(y)}T\chi_{B_L(y)}v\bigr\|\le\varepsilon.
\]
\end{Lem}

\begin{proof}
For each $z\in X$, we choose an ordering of $B_L(z)$ which lists the elements of $B_S(z)$ in the first $\#B_S(z)$ positions.

For each $n\in\{1,\dots,N(L)\}$, let $M_n(\IC)_1$ be the closed unit ball of $M_n(\IC)$ with respect to the operator norm. The set $\prod_{T\in\F}M_n(\IC)_1$ is compact with respect to the maximum norm. Fix a finite partition $\C_n$ of this set into pieces of diameter at most $\varepsilon$. Using the chosen ordering on $B_L(z)$, we view the tuple $(\chi_{B_L(z)}T\chi_{B_L(z)})_{T\in\F}$ as an element of $\prod_{T\in\F}M_n(\IC)_1$ with $n=\#B_L(z)$, and let $G_z\in\C_n$ be the piece containing it. Let $c(z)$ be the triple $\bigl(\#B_L(z),\#B_S(z),G_z\bigr)$. Since $\#B_L(z)\le N(L)$ and each $\C_n$ is finite, $c$ takes only finitely many values, which we enumerate as $1,\dots,\kappa$.

Suppose $c(z)=c(y)$. Then $\#B_L(z)=\#B_L(y)$, $\#B_S(z)=\#B_S(y)$ and $G_z=G_y$. Write $B_L(z)=\{z_1,\dots,z_n\}$ and $B_L(y)=\{y_1,\dots,y_n\}$ in the chosen orderings. We define $\varphi(z_i)=y_i$ so that it takes $B_S(z)$ onto $B_S(y)$ by definition. Then $\langle v^*Tv\delta_{z_k},\delta_{z_i}\rangle=\langle T\delta_{y_k},\delta_{y_i}\rangle$ for all $i,k$, thus the matrix of $v^*\chi_{B_L(y)}T\chi_{B_L(y)}v$ with respect to $(\delta_{z_i})$ is the matrix of $\chi_{B_L(y)}T\chi_{B_L(y)}$ with respect to $(\delta_{y_i})$. Both tuples lie in the same piece $G_z=G_y$ of diameter at most $\varepsilon$, which gives the desired estimate.
\end{proof}

\begin{Lem}\label{lem: depth function}
Let $\P$ be a partition of $X$, let $s\ge0$ and $w>0$, and define $H\colon X\to[0,1]$ by
\[
  H(x)=\min\Bigl\{1,\ \frac{\bigl(d(x,X\setminus\P(x))-s\bigr)_+}{w}\Bigr\},
\]
where $t_+=\max\{t,0\}$ and $d(x,\emptyset)=\infty$. Then
\begin{enumerate}[(1)]
\item $|H(x)-H(y)|\le d(x,y)/w$ for all $x,y\in X$;
\item if $f\colon X\to[0,1]$ satisfies $f\le H$, then $f^{1/2}Tf^{1/2}$ is $\P$-diagonal for every $T\in\B(\ell^2(X))$ with $\prop(T)\le s$;
\item $H(x)=1$ whenever $B_{s+w}(x)\subseteq\P(x)$.
\end{enumerate}
\end{Lem}

\begin{proof}
(1) If $\P(x)=\P(y)$, the claim follows because $x\mapsto d(x,X\setminus\P(x))$ is $1$-Lipschitz on a block and $t\mapsto\min\{1,(t-s)_+/w\}$ is $1/w$-Lipschitz on $\IR$. If $\P(x)\neq\P(y)$, then $y\in X\setminus\P(x)$, so $d(x,X\setminus\P(x))\le d(x,y)$ and hence $0\le H(x)\le d(x,y)/w$. In the same way $0\le H(y)\le d(x,y)/w$, and therefore $|H(x)-H(y)|\le d(x,y)/w$.

(2) The $(x,y)$-entry of $f^{1/2}Tf^{1/2}$ is $f(x)^{1/2}T_{xy}f(y)^{1/2}$. Suppose it is nonzero. Then $H(x)\ge f(x)>0$, so $d(x,X\setminus\P(x))>s$, i.e., $B_s(x)\subseteq\P(x)$. Moreover $T_{xy}\neq0$ gives $d(x,y)\le s$, so $y\in\P(x)$. Hence all nonzero entries of $f^{1/2}Tf^{1/2}$ lie within blocks.

(3) If $B_{s+w}(x)\subseteq\P(x)$, then $d(x,X\setminus\P(x))>s+w$, so $H(x)=1$.
\end{proof}

\subsection{The constructions}\label{sec: construction}

Throughout the rest of this section, we fix $X$ and $\F$ as in Proposition~\ref{pro: twisted partition of unity}, and we choose $s\ge0$ such that $\prop(T)\le s$ for every $T\in\F$. We start with the following constructions:
\begin{enumerate}[(1)]
\item \emph{A slowly varying partition of unity.} Let $\eta>0$ be a small real number to control the error term. Applying Property~A of $X$ (Definition~\ref{def: property A}) with $R=s+1$ and $\varepsilon=\eta$, we obtain $S>0$ and $x\mapsto\xi_x$. For $z\in X$, define $\psi_z\colon X\to[0,1]$ by $\psi_z(x)=\xi_x(z)$. Then $\sum_{z\in X}\psi_z=1$ and $\supp(\psi_z)\subseteq B_S(z)$.

\item \emph{Local types.} Every $T\in\F$ maps $\ell^2(B_S(z))$ into $\ell^2(B_{S+s}(z))$. Let $c$ be the coloring given by Lemma~\ref{lem: local types} for $\F$, $S$, $S+s$ and $\varepsilon=\eta$.

\item \emph{The partition.} Let $w>0$ (it controls the Lipschitz constant of
$H$) and $m\in\IN$. Put $r=S+2s+2w$, and let $\P$ and $\mu$ be given by
Proposition~\ref{pro: geometric matching} for the coloring $c$, $m$ and $r$.
Let $H$ be the depth function of Lemma~\ref{lem: depth function} for $\P$, $s$
and $w$.
\end{enumerate}

With this choice of $r$, we have the following facts.

\begin{Lem}\label{lem: geometric facts}
\begin{enumerate}[(1)]
\item If $z\notin\partial_r\P$, then $(1-H)\psi_z=0$.
\item For every $z\in\partial_r\P$ and $j\in\{1,\dots,m\}$, we have $H=1$ on $B_{S+s+w}(\mu(z,j))$.
\item For distinct $(z,j),(z',j')\in\partial_r\P\times\{1,\dots,m\}$,
\[
  d\bigl(B_{S+s+w}(\mu(z,j)),\,B_{S+s+w}(\mu(z',j'))\bigr)>s.
\]
\end{enumerate}
\end{Lem}

\begin{proof}
(1) Let $z\notin\partial_r\P$, i.e., $B_r(z)\subseteq\P(z)$, and let $x$ satisfy $\psi_z(x)\neq0$. Then $d(x,z)\le S\le r$, so $x\in\P(z)$, and since $S+s+w\le r$,
\[
  B_{s+w}(x)\subseteq B_{S+s+w}(z)\subseteq B_r(z)\subseteq\P(z)=\P(x).
\]
By Lemma~\ref{lem: depth function} (3), $H(x)=1$. Hence $(1-H)\psi_z=0$.

(2) Let $x\in B_{S+s+w}(\mu(z,j))$. By Proposition~\ref{pro: geometric matching} (2), $B_r(\mu(z,j))\subseteq\P(\mu(z,j))$. Since $S+s+w\le r$, this gives $x\in\P(\mu(z,j))$, i.e., $\P(x)=\P(\mu(z,j))$. Since $(S+s+w)+(s+w)=r$,
\[
  B_{s+w}(x)\subseteq B_r(\mu(z,j))\subseteq\P(x),
\]
and Lemma~\ref{lem: depth function} (3) gives $H(x)=1$.

(3) By Proposition~\ref{pro: geometric matching} (3), $d(\mu(z,j),\mu(z',j'))>2r$. Hence any two points of $B_{S+s+w}(\mu(z,j))$ and $B_{S+s+w}(\mu(z',j'))$ are at distance greater than $2r-2(S+s+w)=2s+2w>s$.
\end{proof}

Let $z\in\partial_r\P$ and $j\in\{1,\dots,m\}$. By Proposition~\ref{pro: geometric matching} (1), the points $z$ and $\mu(z,j)$ have the same color. Let $v_{z,j}$ be the partial isometry given by Lemma~\ref{lem: local types} for $z$ and $\mu(z,j)$. Thus $v_{z,j}$ maps $\ell^2(B_{S+s}(z))$ isometrically onto $\ell^2(B_{S+s}(\mu(z,j)))$, maps $\ell^2(B_S(z))$ onto $\ell^2(B_S(\mu(z,j)))$, and satisfies
\[
  \bigl\|\chi_{B_{S+s}(z)}T\chi_{B_{S+s}(z)}-v_{z,j}^*\chi_{B_{S+s}(\mu(z,j))}T\chi_{B_{S+s}(\mu(z,j))}v_{z,j}\bigr\|\le\eta,
\]
for any $T\in\F$. Since $\sup_{z,j}d(z,\mu(z,j))<\infty$, the propagations of the $v_{z,j}$ are uniformly bounded.

For $z\in X$, put $f_z=\bigl((1-H)\psi_z\bigr)^{1/2}$. Then $f_z$ vanishes outside $B_S(z)$, and by Lemma~\ref{lem: geometric facts} (1), $f_z\ne 0$ only if $z\in\partial_r\P$. Hence
\[
  \sum_{z\in\partial_r\P}f_z^2=(1-H)\sum_{z\in X}\psi_z=1-H.
\]
We move the weight $f_z$ to the first partner of $z$ by defining
\[
  K_0=\operatorname{SOT-}\sum_{z\in\partial_r\P}v_{z,1}f_z .
\]
The summands have pairwise orthogonal ranges by Lemma~\ref{lem: geometric facts} (3) and uniformly bounded propagation, so the series converges strongly and $K_0\in C^*_u(X)$.

\begin{Lem}\label{lem: transfer}
\begin{enumerate}[(1)]
\item $K_0^*K_0=1-H$;
\item the range of $K_0$ is contained in
$\ell^2\bigl(\bigsqcup_{z\in\partial_r\P}B_S(\mu(z,1))\bigr)$.
\end{enumerate}
\end{Lem}

\begin{proof}
Let $z\in\partial_r\P$. Since $f_z$ vanishes outside $B_S(z)$ and $v_{z,1}$ maps $\ell^2(B_S(z))$ onto $\ell^2(B_S(\mu(z,1)))$, the range of $v_{z,1}f_z$ lies in $\ell^2(B_S(\mu(z,1)))$, which gives (2). These subspaces are pairwise orthogonal for distinct $z$ by Lemma~\ref{lem: geometric facts} (3), and $v_{z,1}$ is isometric on $\ell^2(B_{S+s}(z))$, which contains the range of $f_z$. Hence, for all $\xi,\xi'\in\ell^2(X)$,
\[
  \langle K_0\xi,K_0\xi'\rangle=\sum_{z\in\partial_r\P}\langle f_z\xi,f_z\xi'\rangle =\langle(1-H)\xi,\xi'\rangle,
\]
which is (1).
\end{proof}

The remaining partners $\mu(z,2),\dots,\mu(z,m)$ carry the constant values $\frac1m,\dots,\frac{m-1}m$. More precisely, we define $g\colon X\to[0,1]$ by
\[
  g(x)=\sum_{z\in\partial_r\P}\ \sum_{j=2}^{m}\frac{j-1}{m} \Bigl(1-\frac{d\bigl(x,B_{S+s}(\mu(z,j))\bigr)}{w}\Bigr)_+ .
\]
The $(z,j)$-summand vanishes outside $B_{S+s+w}(\mu(z,j))$, and these sets are pairwise disjoint by Lemma~\ref{lem: geometric facts} (3). Hence at each point at most one summand is nonzero, and $g\le\frac{m-1}{m}$. Finally, we put
\[
  h=H-g\qquad\text{and}\qquad K=K_0+g^{1/2}.
\]

\begin{Lem}\label{lem: layers}
\begin{enumerate}[(1)]
\item For every $z\in\partial_r\P$ and $j\in\{1,\dots,m\}$, we have $g=\frac{j-1}{m}$ on $B_{S+s}(\mu(z,j))$.
\item $0\le g\le H$ and $0\le h\le H$. In particular, $h^{1/2}Th^{1/2}$ and $g^{1/2}Tg^{1/2}$ are $\P$-diagonal for every $T\in\F$.
\item $K^*K=1-h$ and $KK^*=K_0K_0^*+g$.
\end{enumerate}
\end{Lem}

\begin{proof}
We leave (1) and (2) to the reader. For (3), if $f_z(x)\neq0$, then $H(x)<1$, so by Lemma~\ref{lem: geometric facts} (2) the point $x$ lies outside every $B_{S+s+w}(\mu(z',j))$, and hence $g(x)=0$. Hence $f_zg^{1/2}=0$ for every $z$, and therefore $K_0g^{1/2}=0$. By Lemma~\ref{lem: transfer} (2), the range of $K_0$ lies where $g=0$, so $g^{1/2}K_0=0$. Expanding the products and using
Lemma~\ref{lem: transfer} (1), we get
\[
  K^*K=K_0^*K_0+g=1-H+g=1-h,\qquad KK^*=K_0K_0^*+g. \qedhere
\]
\end{proof}

For $j\in\{1,\dots,m\}$, put
\[
  J_j=\sum_{z\in\partial_r\P}v_{z,j}v_{z,1}^* .
\]
The series converges strongly and $J_j\in C^*_u(X)$. By Lemma~\ref{lem: geometric facts}, the balls $\{B_{S+s}(\mu(z,j))\}_{z\in\partial_r\P, j\in\{1,\dots,m\}}$, are pairwise disjoint. Since $v_{z,j}v_{z,1}^*$ maps $\ell^2(B_{S+s}(\mu(z,1)))$ isometrically onto $\ell^2(B_{S+s}(\mu(z,j)))$ and vanishes elsewhere, we have
\[
  J_j^*J_k=\begin{cases}J_1,& j=k,\\ 0,& j\neq k,\end{cases}\qquad
  J_jJ_j^*=\chi_{\bigcup_{z\in\partial_r\P}B_{S+s}(\mu(z,j))}.
\]
In particular, $J_1$ is the projection onto $\ell^2\bigl(\bigcup_{z\in\partial_r\P}B_{S+s}(\mu(z,1))\bigr)$, and $J_1C^*_u(X)J_1$ is a unital $C^*$-algebra with unit $J_1$. By Lemma~\ref{lem: transfer} (2), $K_0K_0^*$ is a positive contraction in $J_1C^*_u(X)J_1$.

\begin{Lem}\label{lem: rotation}
Let $m\in\IN$. There is a continuous map $u\colon[0,1]\to U(m)$ with $u(0)=1_m$ such that for every $t\in[0,1]$,
\[
  \Bigl\|u(t)^*\diag\bigl(t,\tfrac1m,\dots,\tfrac{m-1}m\bigr)u(t)
  -\diag\bigl(0,\tfrac1m,\dots,\tfrac{m-1}m\bigr)\Bigr\|\le\frac1m .
\]
\end{Lem}

\begin{proof}
For $1\le j\le m-1$ and $\alpha\in\IR$, let $\rho_j(\alpha)\in U(m)$ be the rotation by the angle $\alpha$ in the plane spanned by the $j$-th and $(j+1)$-st coordinates, acting as the identity on the other coordinates. Put $\omega_1=1$ and $\omega_{j+1}=\omega_j\rho_j(\frac\pi2)$ for $1\le j\le m-1$. Since $\rho_j(\frac\pi2)$ exchanges the $j$-th and $(j+1)$-st coordinates up to a sign, conjugating a diagonal matrix by $\omega_j$ moves its first diagonal entry to the $j$-th position and shifts the $2$nd through $j$-th entries one position to the left. In particular,
\[
  \omega_j^*\diag\bigl(t,\tfrac1m,\dots,\tfrac{m-1}m\bigr)\omega_j
  =\diag\bigl(\tfrac1m,\dots,\tfrac{j-1}m,\,t,\,\tfrac jm,\dots,\tfrac{m-1}m\bigr).
\]

We define $u\bigl(\frac{j-1}m\bigr)=\omega_j$ for $1\le j\le m$. On each interval $[\frac{j-1}m,\frac jm]$ with $j\le m-1$, we connect $\omega_j$ to $\omega_{j+1}$ by rotating the $j$-th and $(j+1)$-st coordinates, namely
\[
  u(t)=\omega_j\,\rho_j\bigl(\tfrac\pi2(mt-j+1)\bigr),
\]
and on $[\frac{m-1}m,1]$ we put $u(t)=\omega_m$. Then $u\colon[0,1]\to U(m)$ is continuous and $u(0)=\omega_1=1$.

Let $t\in[\frac{j-1}m,\frac jm]$. Then $u(t)^*\diag(t,\frac1m,\dots,\frac{m-1}m)u(t)$ and $\diag(0,\frac1m,\dots,\frac{m-1}m)$ are block diagonal with respect to the same decomposition of $\{1,\dots,m\}$ into singletons and, if $j\le m-1$, the pair $\{j,j+1\}$. In each pair of corresponding blocks, both blocks are self-adjoint with spectrum in a common interval of length $\frac1m$. Hence their difference has norm at most $\frac1m$.
\end{proof}

Let $u$ be given by Lemma~\ref{lem: rotation}, and write $u(t)=[u_{jk}(t)]_{j,k=1}^m$, so that each $u_{jk}\colon[0,1]\to\IC$ is continuous. Since $K_0K_0^*$ is a positive contraction in the unital $C^*$-algebra $J_1C^*_u(X)J_1$, continuous functional calculus gives elements $u_{jk}(K_0K_0^*)\in J_1C^*_u(X)J_1$. Define
\[
  Q=1-\chi_{\bigcup_{z\in\partial_r\P}\bigcup_{j=1}^mB_{S+s}(\mu(z,j))}
  +\sum_{j,k=1}^mJ_j\,u_{jk}(K_0K_0^*)\,J_k^*
  \qquad\text{and}\qquad V=Q^*K.
\]
One may regard the formula above as a matrix-valued functional calculus. A continuous map $f\colon[0,1]\to M_m(\IC)$ is an element of $M_m(C[0,1])$, and applying functional calculus at $K_0K_0^*$ entrywise
defines a unital $*$-homomorphism
\begin{equation}\label{eq: matrix function calculus}
  M_m(C[0,1])\to M_m\bigl(J_1C^*_u(X)J_1\bigr),\qquad f\mapsto f(K_0K_0^*).
\end{equation}
Here the $m$-matrix coordinates correspond to the $m$ subspaces $\H_j$, where $$\H_j=\ell^2\left(\bigsqcup_{z\in\partial_r\P}B_{S+s}(\mu(z,j))\right),\quad\text{ for any }j=1,\dots,m,$$
By definition, $J_j$ maps $\H_1$ isometrically onto $\H_j$. Accordingly, for $b\in J_1C^*_u(X)J_1$, the operator $J_j\,b\,J_k^*$ maps $\H_k$ into $\H_j$ and vanishes on its orthogonal complement. Hence $\sum_{j,k}J_j\,b_{jk}\,J_k^*$ realizes the matrix $[b_{jk}]$ on $\bigoplus_{j=1}^m\H_j$.

With this picture in mind, we show the following property for $Q$.

\begin{Lem}\label{lem: rotated}
\begin{enumerate}[(1)]
\item $Q$ is a unitary in $C^*_u(X)$;
\item $V^*V=1-h$ and $\|VV^*-g\|\le\frac1m$.
\end{enumerate}
\end{Lem}

\begin{proof}
Let $P$ be the projection onto $\bigoplus_{j=1}^m\H_j$. By the relations for the $J_j$, we have $P=\sum_{j=1}^mJ_jJ_j^*$. As explained above, the map
\[
  \Phi\colon M_m\bigl(J_1C^*_u(X)J_1\bigr)\to PC^*_u(X)P,\qquad
  [b_{jk}]\mapsto\sum_{j,k=1}^mJ_jb_{jk}J_k^*,
\]
is a $*$-isomorphism, with inverse $x\mapsto[J_j^*xJ_k]_{j,k}$, and by definition $Q=(1-P)+\Phi\bigl(u(K_0K_0^*)\bigr)$.

(1) Since $u$ is a unitary in $M_m(C[0,1])$ and $f\mapsto f(K_0K_0^*)$ is a unital $*$-homomorphism, $u(K_0K_0^*)$ is a unitary in $M_m(J_1C^*_u(X)J_1)$. Hence $\Phi\bigl(u(K_0K_0^*)\bigr)$ is a unitary in $PC^*_u(X)P$, and $Q$ is a unitary in $C^*_u(X)$.

(2) By Lemma~\ref{lem: layers} (3), $V^*V=K^*QQ^*K=K^*K=1-h$.

By Lemma~\ref{lem: layers} (1), $g$ equals $\frac{j-1}m$ on $\bigsqcup_{z\in\partial_r\P}B_{S+s}(\mu(z,j))$ for every $j$. Therefore, applying the matrix-valued functional calculus to $K_0K_0^*$ by the continuous functions
\[
  t\mapsto\diag\bigl(t,\tfrac1m,\dots,\tfrac{m-1}m\bigr)
  \qquad\text{and}\qquad
  t\mapsto\diag\bigl(0,\tfrac1m,\dots,\tfrac{m-1}m\bigr),
\]
we get
\[
  K_0K_0^*+Pg=\Phi\Bigl(\diag\bigl(K_0K_0^*,\tfrac1mJ_1,\dots,\tfrac{m-1}mJ_1\bigr)\Bigr),
  \qquad
  Pg=\Phi\Bigl(\diag\bigl(0,\tfrac1mJ_1,\dots,\tfrac{m-1}mJ_1\bigr)\Bigr).
\]
By Lemma~\ref{lem: layers} (3),
\[
  VV^*-g=Q^*(K_0K_0^*+g)Q-g=Q^*(K_0K_0^*+Pg)Q-Pg,
\]
the second equation holds because $g$ commutes with $P$, and $(1-P)g$ is unchanged by conjugation with $Q$, since $Q$ acts as the identity on the range of $1-P$. Consequently, $VV^*-g$ is the image under $\Phi$ of $F(K_0K_0^*)$, where $F\colon[0,1]\to M_m(\IC)$ is the continuous function
\[
  F(t)=u(t)^*\diag\bigl(t,\tfrac1m,\dots,\tfrac{m-1}m\bigr)u(t)
  -\diag\bigl(0,\tfrac1m,\dots,\tfrac{m-1}m\bigr).
\]
Since $\Phi$ is isometric and $*$-homomorphisms are contractive,
\[
  \|VV^*-g\|=\|F(K_0K_0^*)\|\le\sup_{t\in[0,1]}\|F(t)\|\le\frac1m
\]
by Lemma~\ref{lem: rotation}.
\end{proof}

\subsection{Estimates and the proof of Proposition~\ref{pro: twisted partition of unity}}

From the constructions in Section~\ref{sec: construction}, we have already defined the partition $\P$, the function $h, g$ and $V\in C^*_u(X)$ satisfying the first three conditions in Defintion~\ref{def: twisted partition of unity}, see Lemma~\ref{lem: layers} and Lemma~\ref{lem: rotated} . It remains to estimate the commutators in condition (4) of the definition. We first need an elementary estimate.

\begin{Lem}\label{lem: commutator}
Let $s\ge0$, let $(f_z)_{z\in Z}$ be a countable family of bounded real functions on $X$, and let $\varepsilon\ge0$ be such that for any $x,y$ with $d(x,y) \leq s$,
\[
  \sum_{z\in Z}|f_z(x)-f_z(y)|^2\le\varepsilon^2.
\]
Then for every contraction $T$ with $\prop(T)\le s$ and every
$\xi\in\ell^2(X)$,
\[
  \sum_{z\in Z}\|[T,f_z]\xi\|^2\le N(s)^2\varepsilon^2\|\xi\|^2.
\]
In particular, if $f\colon X\to[0,1]$ satisfies $|f(x)-f(y)|\le d(x,y)/w$ for some $w>0$, then $\|[T,f^{1/2}]\|\le N(s)\sqrt{s/w}$.
\end{Lem}

\begin{proof}
We have $([T,f_z]\xi)(x)=\sum_{y\in B_s(x)}T_{xy}\bigl(f_z(y)-f_z(x)\bigr)\xi(y)$. Since $|T_{xy}|\le1$ and $\#B_s(x)\le N(s)$, the Cauchy--Schwarz inequality gives
\[
  \sum_z\|[T,f_z]\xi\|^2
  \le N(s)\sum_{x\in X}\sum_{y\in B_s(x)}\Bigl(\sum_z|f_z(y)-f_z(x)|^2\Bigr)|\xi(y)|^2
  \le N(s)^2\varepsilon^2\|\xi\|^2,
\]
For the last assertion, apply this to the single function $f^{1/2}$ and use $|f(x)^{1/2}-f(y)^{1/2}|^2\le|f(x)-f(y)|\le s/w$ when $d(x,y)\le s$.
\end{proof}

\begin{Lem}\label{lem: first estimates}
For every $T\in\F$:
\begin{enumerate}[(1)]
\item $\|[T,h^{1/2}]\|\le N(s)\sqrt{s/w}$ and
$\|[T,g^{1/2}]\|\le N(s)\sqrt{s/w}$;
\item $\|[T,K_0]\|\le\eta+N(s)\bigl(\sqrt{s/w}+\sqrt\eta\bigr)$.
\end{enumerate}
\end{Lem}

\begin{proof}
(1) Each summand in the definition of $g$ is $1/w$-Lipschitz and nonnegative, and the summands have disjoint supports. Hence $g$ is $1/w$-Lipschitz. Since $H=1$ wherever $g>0$, we have $h=\min\{H,1-g\}$, which is $1/w$-Lipschitz by  Lemma~\ref{lem: depth function} (1). Now apply Lemma~\ref{lem: commutator} we get the desired estimate.

(2) Fix $z\in\partial_r\P$. Since $\prop(T)\le s$, the operators $Tf_z$ and $f_zT$ have range in $\ell^2(B_{S+s}(z))$, and $Tv_{z,1}f_z$ has range in $\ell^2(B_{S+s}(\mu(z,1)))$. Using $v_{z,1}v_{z,1}^*=\chi_{B_{S+s}(\mu(z,1))}$, one checks that
\[
  [T,v_{z,1}f_z]
  =v_{z,1}\bigl(v_{z,1}^*\chi_{B_{S+s}(\mu(z,1))}T\chi_{B_{S+s}(\mu(z,1))}v_{z,1}
  -\chi_{B_{S+s}(z)}T\chi_{B_{S+s}(z)}\bigr)f_z+v_{z,1}[T,f_z].
\]
Summing over $z$, we write $[T,K_0]=A_1+A_2$, where $A_1$ and $A_2$ are the sums of the first and the second terms. In both sums, the $z$-th term has range in $\ell^2(B_{S+s}(\mu(z,1)))$, and these subspaces are pairwise orthogonal. By Lemma~\ref{lem: local types} and Lemma~\ref{lem: transfer}, $\sum_z\|f_z\xi\|^2\le\|\xi\|^2$ thus we get $\|A_1\|\le\eta$.

For $A_2$, let $d(x,y)\le s$. Writing $f_z=(1-H)^{1/2}\psi_z^{1/2}$ and using the triangle inequality in $\ell^2$ together with $\sum_z\psi_z=1$, we get
\[
  \Bigl(\sum_z|f_z(x)-f_z(y)|^2\Bigr)^{1/2}
  \le\bigl|(1-H(x))^{1/2}-(1-H(y))^{1/2}\bigr|
  +\Bigl(\sum_z\bigl|\psi_z(x)^{1/2}-\psi_z(y)^{1/2}\bigr|^2\Bigr)^{1/2}.
\]
Since $|\alpha^{1/2}-\beta^{1/2}|^2\le|\alpha-\beta|$ for $\alpha,\beta\ge0$, the first term is at most $|H(x)-H(y)|^{1/2}\le\sqrt{s/w}$, and the second is at most $\bigl(\sum_z|\xi_x(z)-\xi_y(z)|\bigr)^{1/2}<\sqrt\eta$ since $d(x,y)<s+1$. Since $v_{z,1}$ is isometric on $\ell^2(B_{S+s}(z))$, Lemma~\ref{lem: commutator} gives
\[
  \|A_2\xi\|^2=\sum_z\|[T,f_z]\xi\|^2
  \le N(s)^2\bigl(\sqrt{s/w}+\sqrt\eta\bigr)^2\|\xi\|^2 .
\]
Adding the two bounds gives (2).
\end{proof}

To estimate the commutators with $Q$, we need the continuity of the matrix-valued functional calculus with respect to commutators.

\begin{Lem}\label{lem: continuity}
Let $m\in\IN$ and let $f\colon[0,1]\to M_m(\IC)$ be continuous. For every $\varepsilon>0$, there exists $\delta>0$ such that for every unital $C^*$-algebra $\mathcal D$, every positive contraction $e\in\mathcal D$ and every contraction $a\in\mathcal D$ with $\|[e,a]\|<\delta$, we have
\[
  \bigl\|[f(e),\diag(a,\dots,a)]\bigr\|<\varepsilon\qquad\text{in }M_m(\mathcal D),
\]
where $f(e)$ is the matrix-value functional calculus as in \eqref{eq: matrix function calculus}.
\end{Lem}

\begin{proof}
By the Weierstrass theorem, there are polynomials $p_{jk}(t)=\sum_{n\ge0}c_{jk,n}t^n$ with $\sup_{t\in[0,1]}|f_{jk}(t)-p_{jk}(t)|<\frac\varepsilon{4m}$. Put $M=\max_{j,k}\sum_{n\ge1}n|c_{jk,n}|$ and $\delta=\frac\varepsilon{2m(M+1)}$. Since $[e^n,a]=\sum_{i=0}^{n-1}e^i[e,a]e^{n-1-i}$, we have $\|[e^n,a]\|\le n\|[e,a]\|$, and hence
\[
  \|[f_{jk}(e),a]\|\le\frac{2\varepsilon}{4m}+M\|[e,a]\|<\frac\varepsilon m .
\]
The $(j,k)$-entry of $[f(e),\diag(a,\dots,a)]$ is $[f_{jk}(e),a]$, and the norm of an $m\times m$ operator matrix is at most $m$ times the largest norm of its entries.
\end{proof}

\begin{Lem}\label{lem: Q estimate}
For every $\varepsilon>0$, there exists $\delta>0$ such that if $\|[T,K_0]\|<\delta$ for every $T\in\F$, then
\[
  \|[Q,T]\|<\varepsilon+4\eta\qquad\text{for every }T\in\F.
\]
\end{Lem}

\begin{proof}
Let $\delta‘$ be given by Lemma~\ref{lem: continuity} for $f=u$ and $\varepsilon$, and fix $T\in\F$. Put $\delta=\delta'/2$.

\emph{Step 1: Localizing $Q-1$.} Since $\sum_jJ_jJ_1J_j^*=\sum_jJ_jJ_j^*$, we have
\[
  Q-1=\sum_{j,k=1}^mJ_j\bigl(u_{jk}(K_0K_0^*)-\delta_{jk}J_1\bigr)J_k^*.
\]
Each function $u_{jk}-\delta_{jk}$ vanishes at $0$ because $u(0)=\mathbb{1}$. By the Weierstrass approximation theorem, it is a uniform limit on $[0,1]$ of polynomials without constant term. Hence $u_{jk}(K_0K_0^*)-\delta_{jk}J_1$ is a norm limit of polynomials in $K_0K_0^*$ without constant term. By Lemma~\ref{lem: transfer} (2), these operators vanish on the orthogonal complement of $\ell^2\bigl(\bigsqcup_zB_S(\mu(z,1))\bigr)$ and have range in it. Since $J_j$ maps $\ell^2(B_S(\mu(z,1)))$ onto $\ell^2(B_S(\mu(z,j)))$, it follows that $Q-1$ vanishes on the orthogonal complement of $\ell^2\bigl(\bigsqcup_{z,j}B_S(\mu(z,j))\bigr)$ and has range in it.

\emph{Step 2: Reduction to the layers.} Let $P$ be the projection onto $\bigoplus_{j=1}^m\H_j$. Since $\prop(T)\le s$, the operator $T$ maps $\ell^2(B_S(\mu(z,j)))$ into $\ell^2(B_{S+s}(\mu(z,j)))\subseteq P\ell^2(X)$, and so does $T^*$. By Step 1, $(1-P)T(Q-1)=0$ and $(Q-1)T(1-P)=0$. By the definition of $Q$, we have $Q-1=P(Q-1)P$, hence
\[
  [Q,T]=[Q-1,T]=[Q-1,PTP].
\]
By Lemma~\ref{lem: geometric facts} (3), $J_j^*TJ_k=0$ for $j\neq k$, so $PTP=\Phi\bigl(\diag(J_1^*TJ_1,\dots,J_m^*TJ_m)\bigr)$. Since $Q-1=\Phi\bigl(u(K_0K_0^*)\bigr)-P$, we obtain
\[
  \|[Q,T]\|=\bigl\|[u(K_0K_0^*),\diag(J_1^*TJ_1,\dots,J_m^*TJ_m)]\bigr\|.
\]

\emph{Step 3: Comparing the layers.} For every $j$, the operator $J_j^*TJ_j$ is block diagonal with respect to the subspaces $\ell^2(B_{S+s}(\mu(z,1)))$, $z\in\partial_r\P$, with blocks $v_{z,1}\bigl(v_{z,j}^*\chi_{B_{S+s}(\mu(z,j))}T\chi_{B_{S+s}(\mu(z,j))}v_{z,j}\bigr)v_{z,1}^*$. By Lemma~\ref{lem: local types}, each block lies within $\eta$ of $v_{z,1}\chi_{B_{S+s}(z)}T\chi_{B_{S+s}(z)}v_{z,1}^*$. Hence $\|J_j^*TJ_j-J_1TJ_1\|\le2\eta$, and
\[
  \|[Q,T]\|\le\bigl\|[u(K_0K_0^*),\diag(J_1TJ_1,\dots,J_1TJ_1)]\bigr\|+4\eta .
\]

Since $K_0K_0^*=J_1K_0K_0^*J_1$ and $\F^*=\F$,
\[
  \|[K_0K_0^*,J_1TJ_1]\|\le\|[K_0K_0^*,T]\|
  \le\|[K_0,T]\|+\|[K_0,T^*]\|<2\delta .
\]
Applying Lemma~\ref{lem: continuity} in the unital $C^*$-algebra $J_1C^*_u(X)J_1$ to $e=K_0K_0^*$ and $a=J_1TJ_1$, which satisfy $\|[e,a]\|<2\delta=\delta'$, gives the claim.
\end{proof}

\begin{proof}[Proof of Proposition~\ref{pro: twisted partition of unity}]
Let $\varepsilon>0$. Choose an integer $m$ with $\frac1m\le\varepsilon$, let $\delta$ be given by Lemma~\ref{lem: Q estimate} for $\frac\varepsilon2$ (this depends only on $m$ and $\varepsilon$), and choose $\eta>0$ and $w>0$ such that
\[
  \beta = \eta+N(s)\bigl(\sqrt{s/w}+\sqrt\eta\bigr)<\min\Bigl\{\delta,\frac\varepsilon{12}\Bigr\}.
\]
Carry out the construction of Section~\ref{sec: construction} with these $\eta$, $w$ and $m$. We verify that $h$, $g$ and $V$ form a twisted partition of unity subordinate to $\P$ for $\F$ and $\varepsilon$. Conditions (1) and (2) follow from Lemma~\ref{lem: rotated}, since $\frac1m\le\varepsilon$, and condition (3) is Lemma~\ref{lem: layers} (2).

For condition (4), let $T\in\F$. By Lemma~\ref{lem: first estimates}, $\|[T,h^{1/2}]\|\le\beta<\varepsilon$, $\|[T,K_0]\|\le\beta<\delta$, and $\|[T,K]\|\le\|[T,K_0]\|+\|[T,g^{1/2}]\|\le2\beta$. By Lemma~\ref{lem: Q estimate} applied to $T^*\in\F$, $\|[Q^*,T]\|=\|[Q,T^*]\|<\frac\varepsilon2+4\eta\le\frac\varepsilon2+4\beta$. Since $\|K\|\le1$ by Lemma~\ref{lem: layers} (3),
\[
  \|[V,T]\|\le\|[K,T]\|+\|[Q^*,T]\|<6\beta+\frac\varepsilon2<\varepsilon. \qedhere
\]
\end{proof}

\section{Proof of the main theorem}\label{sec: proof of the main theorem}

Recall \cite[Definition~2.1]{WZ10} that a $C^*$-algebra $A$ has \emph{nuclear dimension at most $n$} if for every finite set $\F\subseteq A$ and every $\varepsilon>0$ there are finite-dimensional $C^*$-algebras $E_0,\dots,E_n$, completely positive contractions $\psi_i\colon A\to E_i$ and completely positive contractive order zero maps $\beta_i\colon E_i\to A$ such that $\|\sum_{i=0}^n\beta_i\psi_i(T)-T\|<\varepsilon$ for every $T\in\F$. Here a completely positive map has \emph{order zero} if it maps orthogonal positive elements to orthogonal elements. Every $*$-homomorphism has order zero, and the composition of an order zero map with a $*$-homomorphism has order zero.

Let $\P$ be a partition of $X$ and $Y\subseteq X$. The $\P$-diagonal operators supported on $\ell^2(Y)$ form a $C^*$-algebra isomorphic to $\prod_{D\in\P,\,D\cap Y\neq\emptyset}\B(\ell^2(D\cap Y))$. Since $\#(D\cap Y)\le\sup_{D\in\P}\#D<\infty$, this algebra is AF by \cite[Lemma~8.4]{WZ10}, and hence it has nuclear dimension zero \cite[Remarks~2.2]{WZ10}.

Now, we are ready to prove our main result.

\begin{proof}[Proof of Theorem~\ref{thm: main theorem}]
Let $\F\subseteq C^*_u(X)$ be finite and let $\varepsilon>0$. We may assume that $\F$ consists of contractions. For each $T\in\F$, choose a contraction $T'\in C^*_u(X)$ of finite propagation with $\|T-T'\|<\varepsilon/8$, and let $\F'$ consist of all $T'$ and their adjoints. Choose $\varepsilon'>0$ with $2\varepsilon'+6\sqrt{\varepsilon'}<\varepsilon/2$. By Proposition~\ref{pro: twisted partition of unity}, there is a twisted partition of unity $h,g,V$ subordinate to some partition $\P$ for $\F'$ and $\varepsilon'$.

Let $B_0$ be the $C^*$-algebra of all $\P$-diagonal operators on $\ell^2(X)$, which is AF and a subalgebra of $C^*_u(X)$. Let $\psi_0(T)=\sum_{D\in\P}\chi_Dh^{1/2}Th^{1/2}\chi_D$, and let $\beta_0$ be the inclusion of $B_0$ into $C^*_u(X)$. For $T'\in\F'$, the operator $h^{1/2}T'h^{1/2}$ is $\P$-diagonal, so $\beta_0\psi_0(T')=h^{1/2}T'h^{1/2}$ and
\[
  \|\beta_0\psi_0(T')-hT'\|=\|h^{1/2}[T',h^{1/2}]\|<\varepsilon'.
\]

A naive choice for the second map would be $T\mapsto V^*\bigl(\sum_{D\in\P}\chi_Dg^{1/2}Tg^{1/2}\chi_D\bigr)V$. This has two defects. First, $x\mapsto V^*xV$ is in general not of order zero, since $V$ is not a partial isometry. Second, the weight is counted twice. Since $VV^*\approx g$, we have $V^*g^{1/2}Tg^{1/2}V\approx(1-h)^2T$ instead of $(1-h)T$. To solve this, we need to find a partial isometry $W$ such that $V\approx g^{1/2}W$. Then $x\mapsto W^*xW$ is a $*$-homomorphism, and $W^*g^{1/2}Tg^{1/2}W\approx V^*TV\approx(1-h)T$. We obtain $W$ from the polar decomposition of $V$. To make sure that $W$ lies in $C^*_u(X)$, we first cut off the region where $g$ is small.

Let $p$ be the characteristic function of $\{x\in X\mid g(x)\ge2\varepsilon'\}$. Since $\|VV^*-g\|\le\varepsilon'$ and $gp\ge2\varepsilon'p$, we have $pVV^*p\ge\varepsilon'p$, so $pVV^*p$ is invertible in $pC^*_u(X)p$. Put $W=(pVV^*p)^{-1/2}pV$, with the inverse taken in $pC^*_u(X)p$. Then $WW^*=p$. In particular $pW=W$, and $pV=(pVV^*p)^{1/2}W$. Since $pW=W$ and $g$ commutes with $p$, we have
$g^{1/2}W-pV=\bigl((gp)^{1/2}-(pVV^*p)^{1/2}\bigr)W$. Using
$\|\alpha^{1/2}-\beta^{1/2}\|\le\|\alpha-\beta\|^{1/2}$ for positive
$\alpha,\beta$ and $gp-pVV^*p=p(g-VV^*)p$, we get
\[
  \|g^{1/2}W-pV\|\le\|p(g-VV^*)p\|^{1/2}\le\sqrt{\varepsilon'}.
\]
Moreover, since $g<2\varepsilon'$ on the support of $1-p$,
\[
  \|(1-p)V\|^2=\|(1-p)VV^*(1-p)\|
  \le\|(1-p)g(1-p)\|+\|VV^*-g\|\le3\varepsilon'.
\]
Writing $V=pV+(1-p)V$, we obtain
$\|g^{1/2}W-V\|\le\sqrt{\varepsilon'}+\sqrt{3\varepsilon'}<3\sqrt{\varepsilon'}$. The corner $pB_0p$ is AF by the discussion above. Let $\psi_1(T)=\sum_{D\in\P}\chi_D\,p\,g^{1/2}Tg^{1/2}p\,\chi_D$ and $\beta_1(x)=W^*xW$. Since $WW^*=p$, $\beta_1$ is a $*$-homomorphism on $pB_0p$. For $T'\in\F'$, the operator $p\,g^{1/2}T'g^{1/2}p$ is $\P$-diagonal, so $\beta_1\psi_1(T')=(g^{1/2}W)^*T'(g^{1/2}W)$. Since $(1-h)T'=V^*VT'$ and $\|V\|,\|g^{1/2}W\|\le1$,
\[
  \|\beta_1\psi_1(T')-(1-h)T'\|
  \le\|(g^{1/2}W)^*T'(g^{1/2}W)-V^*T'V\|+\|V^*[T',V]\|
  <6\sqrt{\varepsilon'}+\varepsilon'.
\]

The maps $\psi_i$ are completely positive contractions and the maps $\beta_i$ are $*$-homomorphisms. Adding the two colors, we get $\|\beta_0\psi_0(T')+\beta_1\psi_1(T')-T'\|<2\varepsilon'+6\sqrt{\varepsilon'}<\varepsilon/2$ for $T'\in\F'$. Since $G=\beta_0\psi_0+\beta_1\psi_1$ has norm at most $2$, for $T\in\F$ we obtain
\[
  \|G(T)-T\|\le\|G\|\,\|T-T'\|+\|G(T')-T'\|+\|T'-T\|
  <\frac{2\varepsilon}8+\frac\varepsilon2+\frac\varepsilon8<\varepsilon .
\]
Since finite-dimensional $C^*$-algebras are injective and the restriction of an order zero map to a subalgebra has order zero, it follows from the definition that we may replace the finite-dimensional $C^*$-algebras $E_i$ by AF algebras.
\end{proof}

\begin{proof}[Proof of Corollary~\ref{cor: trichotomy}]
Suppose first that $X$ has Property~A. By Theorem~\ref{thm: main theorem},
$\dimnuc C^*_u(X)\le1$, and by \cite[Corollary~1.5]{LW18},
$\dimnuc C^*_u(X)=0$ if and only if $\asdim X=0$. Note that $\asdim X=0$
implies Property~A, so the first two cases are covered.

Suppose now that $X$ does not have Property~A. Then $C^*_u(X)$ is not nuclear
\cite{STY02}, \cite{BO08}. Since finite nuclear dimension implies nuclearity
\cite[Remarks~2.2]{WZ10}, we get $\dimnuc C^*_u(X)=\infty$.
\end{proof}

\begin{Rem}\label{rem: unbounded geometry}
The bounded geometry assumption in Corollary~\ref{cor: trichotomy} is harmless. Let $X$ be a discrete metric space which does not have bounded geometry, and let $C^*_u(X)$ be the norm closure of the bounded operators on $\ell^2(X)$ of finite propagation. Then $\dimnuc C^*_u(X)=\infty$. Indeed, choose $r>0$ such that $\sup_{x\in X}\#B_r(x)=\infty$. Inductively, we find pairwise disjoint finite sets $S_1,S_2,\dots$ with $\#S_n=n$, each contained in a ball of radius $r$. The operators which are block diagonal with respect to the sets $S_n$ and vanish on $\ell^2\bigl(X\setminus\bigsqcup_nS_n\bigr)$ have propagation at most $2r$. Hence $C^*_u(X)$ contains a copy of $\prod_{n\ge1}M_n(\IC)$, which is not exact \cite{BO08}. Since exactness passes to $C^*$-subalgebras and finite nuclear dimension implies nuclearity, we get $\dimnuc C^*_u(X)=\infty$.
\end{Rem}

By Corollary~\ref{cor: trichotomy} and Remark~\ref{rem: unbounded geometry}, the nuclear dimension of a uniform Roe algebra can only take the values in $\{0, 1, \infty\}$. The same trichotomy holds for simple $C^*$-algebras \cite[Corollart C]{CETWW21}, \cite[Corollary C]{CE20}. These two classes of $C^*$-algebras are of very different nature, and it is natural to ask what they have in common.

\begin{Que}
For which class of $C^*$-algebras $A$ does one have
$\dimnuc A\in\{0,1,\infty\}$?
\end{Que}

\bibliographystyle{myalpha}
\bibliography{ref}

\end{document}